\documentclass[preprint,sort&compress,a4paper]{elsarticle}
\usepackage{graphicx}
\usepackage[latin1]{inputenc}

\usepackage{array}
\usepackage{xcolor}
\usepackage{amsmath,amsfonts,latexsym,amssymb,amsbsy}

\usepackage{subcaption}
\usepackage{graphicx}

\usepackage{float}
\usepackage{hyperref}

\usepackage{framed}

		\usepackage{pgfplots}
		\usetikzlibrary{plotmarks}
		\usetikzlibrary{external}
		\pgfplotsset{compat=1.3}
		\newlength\figurewidth 
		\newlength\figureheight
	\definecolor{mycolor1}{rgb}{0.549019634723663,0.141176477074623,0.850980401039124}%
\definecolor{mycolor2}{rgb}{0.622055131485066,0.350952380892271,0.513249539867053}%
\definecolor{mycolor3}{rgb}{0,0.498039215803146,0}%
\definecolor{mycolor4}{rgb}{0.337254911661148,0.898039221763611,0.368627458810806}%
\definecolor{mycolor5}{rgb}{1,0,1}%
\definecolor{mycolor6}{rgb}{0.392227019534168,0.655477890177557,0.171186687811562}%

\pgfplotsset{
	styleRKexplicit/.style=
	{color=mycolor1,solid,mark=o,mark options={solid},line width=1pt},
	styleRKGL/.style=
	{color=mycolor2,dotted,mark=+,mark size=3pt,mark options={solid},line width=1pt},
	styleRKN/.style=
	{color=black,solid,mark=square*,line width=1pt},
  stylePhiFiveSix/.style=
	{color=blue,dashed,mark=x,mark size=3pt,mark options={solid},line width=1pt},
 stylePhiOneSix/.style=
	{color=mycolor3,dash pattern=on 1pt off 3pt on 3pt off 3pt,mark=diamond,mark size=3pt,mark options={solid},line width=1pt},
 stylePhiTwoSix/.style=
	{color=red,dashed,line width=1pt},
	stylePhiThreeSix/.style=
	{color=mycolor4,solid,mark=*,mark options={solid},line width=1pt},
	stylePhiFiveEight/.style=
	{color=mycolor5,solid,mark=triangle*,mark size=3pt}
	}
	
\tikzset{external/force remake=false} 

\newcommand{\CC}{\mathcal{C}}
\newcommand{\cost}{cost}

\begin{document}

\begin{frontmatter}

\title{Symplectic integrators for the matrix Hill's equation and its applications to engineering models}
\journal{}

\author[lat]{Philipp Bader}\corref{cor1}\ead{p.bader@latrobe.edu.au}
\author[imm]{Sergio Blanes}\ead{serblaza@imm.upv.es}
\author[imm]{Enrique Ponsoda}\ead{eponsoda@imm.upv.es}
\author[zaza]{Muaz Seydao\u{g}lu}\ead{muasey@upvnet.upv.es}
\address[lat]{Department of Mathematics and Statistics, La Trobe University, 3086 Bundoora VIC, Australia}
\address[imm]{Instituto de Matem\'{a}tica Multidisciplinar, 
Universitat Polit\`{e}cnica de Val\`{e}ncia. Spain}
\address[zaza]{Department of Mathematics, Faculty of Art and Science, Mu\c{s} Alparslan University, 49100, Mu\c{s}, Turkey}
\cortext[cor1]{Corresponding author}

\begin{abstract}
We consider the numerical integration of the matrix Hill's equation. Parametric resonances can appear and this property is of great interest in many different physical applications. Usually, the Hill's equations originate from a Hamiltonian function and the fundamental matrix solution is  a symplectic matrix. This is a very important property to be preserved by the numerical integrators. In this work we present new sixth-and eighth-order symplectic exponential integrators that are tailored to the Hill's equation. The methods are based on an efficient symplectic approximation to the exponential of high dimensional coupled autonomous harmonic oscillators and yield accurate results for oscillatory problems at a low computational cost. Several numerical examples illustrate the performance of the new methods.

\end{abstract}
\begin{keyword}
matrix Hill's equation; Mathieu equation; parametric resonance; symplectic integrators; Magnus expansion
\MSC 65L07; 65L05  ; 65Z05 
\end{keyword}
\end{frontmatter}

\section{Introduction}
The study of the potential of a charged particle moving in the electric field of a quadrupole, without considering the effects of the induced magnetic field,  leads to the equations of motion 
\begin{displaymath}
\frac{d^2}{dt^2}
\begin{pmatrix}
x_{1} \\ x_{2}  \\ x_{3} 
\end{pmatrix}
=
\begin{pmatrix}
\displaystyle  - \frac{2e}{m d} \left( E_{0} + E_{1} \cos (t) \right) & 0 & 0 \\
0 & \displaystyle  \frac{2e}{m d} \left( E_{0} + E_{1} \cos (t) \right)    & 0 \\
0 & 0 & 0
\end{pmatrix}
\begin{pmatrix}
x_{1}  \\ x_{2}  \\ x_{3} 
   \end{pmatrix},
\end{displaymath}
where $e$ is the charge of the particle, $m$ is the mass and $d$ is the minimum distance from the electrode to the $x_3$-axis
(the direction in which the particle is traveling), and $ \, E(t) = E_{0} + E_{1} \cos (t) \, $ is the electric field which is radio-frequency-modulated, see \cite{ruby}.

The decoupled equations for $x_{i}$, $ i = 1, 2$, are the well-known Mathieu equations
\begin{displaymath}
x'' (t) + \left( \alpha - \beta \cos (t) \right) x (t) = 0 \, ;
\end{displaymath}
where $ \alpha $ and $ \beta $ are constant parameters and each prime denotes a derivative in time.

Stable trajectories are guaranteed only in certain regions in the $\alpha-\beta$ plane and hence, the the motion of the particles can be regulated by $ E_{0} $ and $ E_{1} $ by placing them within or outside of said stability regions.

This property allows to filter particles and is the underlying principle of quadrupole mass spectrometry. The setup is also known as Paul's trap \cite{paul90etf}. This important application motivated a great interest in the study of periodically time-dependent systems and the effect of parametric resonances.

In this work, we consider the more general problem which is a matrix version of the so called Hill's equation
\begin{equation}\label{hill1}
x'' + M(t) x = 0 \, ,
\end{equation}
where $  t \in \mathbb{R} $, $ x \in \mathbb{C}^{r} $ and $ M(t) $ is a periodic $ r \times r $ matrix valued function with period $T$. Frequently, eq. (\ref{hill1}) is written as
\begin{equation}\label{hill2}
x'' + \left(A+\sum_{n=1}^{\infty} B_n\cos(nt)\right) x = 0 \, ,
\end{equation}
which reduces - when truncating after the first term - to the Mathieu equation.

The Hill's equation has many applications in practical periodically variable systems like the study of  qua\-dru\-po\-le mass filters and qua\-dru\-po\-le devices \cite{drewsen00hlp,paul90etf}, microelectromechanical systems \cite{turner98fpr}, parametric resonances in Bose-Einstein condensates \cite{cairncross14pri,garcia99epr}, spatially linear electric fields, dynamic buckling of structures, electrons in crystal lattices, waves in periodic media, etc. (see also \cite{magnus66he,major05cpt,mclachlan65taa,richards,werth09cpt} and references therein).

In \cite{inicio}, the theory of Floquet is applied in order to study equation (\ref{hill1}). 

The majority of published works discusses analytical methods to find the stability regions. However, the computation of the trajectories is of great interest in many cases, and this has to be carried out numerically. 
In most cases, Hill's equation originates from a Hamiltonian system where $M^T=M$, and the fundamental matrix solution is a symplectic matrix. This property plays a fundamental role on the stability of the system.

For Hamiltonian systems it is essential to preserve the symplectic structure of the exact solution for both the study of stability regions and for very long time numerical integrations. Therefore, our goal in this work is to design symplectic methods based on Magnus expansions which will efficiently integrate Hill's equation. We analyze new methods that are closely related to commutator-free methods \cite{alvermann11hoc,blanes02psp} and integrators for $N$th-order linear systems \cite{bbcp} but they are tailored for solving the relevant matrix Hill's equation. 
The first step in this undertaking is to write Hill's eq. (\ref{hill1}) as a first-order system
\begin{equation}\label{MatrixHill}
\frac{d}{dt}\begin{pmatrix} x \\ x' \end{pmatrix}
 =  \begin{pmatrix}
0 & I \\ - M(t) & 0
\end{pmatrix} \begin{pmatrix} x \\ x'  \end{pmatrix},
\end{equation}
and we propose new sixth- and eighth-order symplectic methods for its solution. The most efficient sixth-order symplectic method that we have found is summarized in Table~\ref{method62} which describes the evolution operator $\Psi_n$ for a step of size $h$ from $t_n$ to $t_{n+1}=t_n+h$
\begin{equation*}\label{MapHill}
\begin{pmatrix} x \\ x' \end{pmatrix}_{n+1}
 =  \Psi_n 
\begin{pmatrix} x \\ x'  \end{pmatrix}_n.
\end{equation*}

\begin{table}
\caption{Two-exponential sixth-order symplectic method using the Gauss-Legendre quadrature rule for solving the Hill's equation (\ref{MatrixHill}).}
\label{method62}
 \begin{framed}
$c_1=\frac{5-\sqrt{15}}{10}, \ c_1=\frac12, \ c_3=\frac{5+\sqrt{15}}{10}$
\[
  M_1=M(t_n+c_1h), \qquad  M_2=M(t_n+c_2h), \qquad  M_3=M(t_n+c_3h), 
\]
\[
 K=M_1-M_3, \qquad  L=-M_1+2M_2-M_3, \qquad  F=K^2,
\]
\[
\begin{array}{ll}
	C_1=-\frac{\sqrt{15}}{180}K + \frac1{18}L+\frac1{12960}F , \qquad
	& 
	D_1=-M_2 -\frac{4}{3\sqrt{15}}K + \frac1{9}L \\
  C_2=+\frac{\sqrt{15}}{180}K + \frac1{18}L+\frac1{12960}F,
	&  
	D_2=-M_2 +\frac{4}{3\sqrt{15}}K + \frac1{9}L
\end{array}
\]
\begin{equation*}
\Psi_n =\begin{pmatrix} I & 0 \\ hC_{2} & I \end{pmatrix}
\exp \left(\frac{h}{2} \begin{pmatrix} 0 & I \\ D_{2} & 0 \end{pmatrix}\right)
\exp \left(\frac{h}{2} \begin{pmatrix} 0 & I \\ D_{1} & 0 \end{pmatrix}\right)
\begin{pmatrix} I & 0 \\ hC_{1} & I \end{pmatrix}
\end{equation*}
\end{framed}
\end{table}
Notice that $M(t)$ is evaluated at the nodes of the sixth-order Gauss-Legendre quadrature rule, but the method can easily be adapted to any other quadrature rule.

Furthermore, the matrices $C_i,D_i, \ i=1,2$ in the algorithm require only three evaluations of $M(t)$, several linear combinations of them and one matrix-matrix product.
Each exponential can be computed up to an accuracy of order $h^{2m}$ using $m-1$ matrix-matrix products and symplecticity can be preserved at an equivalent extra cost of $2+\frac13$ products (see the Appendix).
This method provides very high accuracy for oscillatory problems while requiring much less evaluations of $M(t)$ per step at a slightly higher computational cost per step when compared to the most efficient explicit symplectic Runge-Kutta-Nystr\"om methods from the literature.
Moreover, the method provides the exact solution in the autonomous case.

In this work, we analyze sixth-order methods with one to three exponentials and eighth-order methods with four and five exponentials. 

We also remark that for the non-homogeneous linear problem
\begin{equation*}\label{HillNonHom}
x''(t) +M(t) x(t) = f(t),
\end{equation*}
the stability of the system is independent of the non-homogeneous term $f(t)$ (although the solution can strongly depend on it). This equation is frequently solved by variation of constants, but a more efficient procedure is to write the system as an homogeneous one by extending it as follows \cite{blapon}
\begin{equation*}\label{PertMatrixHill}
\frac{d}{dt}\begin{pmatrix} x \\ x' \\ 1 \end{pmatrix}
 =  \begin{pmatrix} 
0_{r \times r} & I_{r \times r} & 0_{r \times 1}  \\
- M(t) & 0_{r \times r} & f(t) \\
 0_{1 \times r}  & 0_{1 \times r} & 0_{1 \times 1} 
\end{pmatrix}
\begin{pmatrix} x \\ x'  \\ 1  \end{pmatrix}.
\end{equation*}
Now it is obvious that the eigenvalues of the fundamental matrix solution will not depend on $f(t)$. 
The methods presented in this work can be applied to this system with minor changes and the computational cost remains essentially the same since it increases only by a low number of vector-matrix multiplications due to the non-homogeneous term. 
The method allows for a different treatment of the matrix $M(t)$ and the vector $f(t)$ \cite{blapon} if required (they can even be evaluated at different nodes).

In section~2, we discuss the numerical integration of Hill's equation by standard methods and briefly introduce the background necessary for the derivation of new methods in section~3.
In section~4, we evaluate the performance of the new methods using numerical examples.

\section{Numerical integration for one period}

Let us consider eq. \eqref{hill1} and the equivalent first order system
\begin{equation}\label{firstordersystem}
z'(t) = A(t) z(t) \, , \ \ \
A(t) =  \begin{pmatrix}
0 & I \\
- M(t) & 0
\end{pmatrix} ,
\end{equation}
with $ \, z(0) =z_{0}\in\mathbb{C}^{2r}$ and periodic  $A(t+T)=A(t)$.
Let $ \Phi(t) $ denote the fundamental matrix solution of (\ref{firstordersystem}), i.e., ${z(t) = \Phi (t) z(0)}$,  then
\begin{equation}\label{Fundamental}
\Phi' (t) = A(t) \Phi (t) , \quad \Phi(0) = I_{2r \times 2r} .
\end{equation}
Floquet theory tells us that  $\Phi (t+T)=\Phi (t)\Phi (T)$. Hence, 
\begin{displaymath}
z(T) = \Phi (T) z_{0} \, , \ \
z(2T) = \Phi (T) z(T) = \Phi^{2} (T) z_{0} \, , \ \ldots \, , \,
z(nT) = \Phi^{n} (T) z_{0} \, ,
\end{displaymath}
and the system is stable if the eigenvalues $\{\lambda_1,\ldots,\lambda_{2r}\}$ of $\Phi (T)$ lie in the unit disk, i.e., $|\lambda_i|\leq 1, \ i=1,\ldots,2r$. Notice that for $|\lambda_i|<1$, the system is asymptotically stable in the direction of the eigenvector associated to this eigenvalue.

Most systems are Hamiltonian (then $M^T=M$)  and $A(t)$ defined by \eqref{firstordersystem} belongs to the symplectic Lie algebra, i.e., $A(t)^T\,J+J\,A(t)=0, \ \forall t$ where
\[
  J =  \begin{pmatrix}
0 & I \\
- I & 0
\end{pmatrix}
\]
is the fundamental symplectic matrix. Furthermore, $ \Phi(t)$ is a symplectic matrix, i.e., $\Phi(t)^T\,J\,\Phi(t)=J,  \forall t$, and $\det \Phi(t)=1$. 
The eigenvalues of $\Phi(T)$ occur in reciprocal pairs,  $\{\lambda,\lambda^*, 1/\lambda, 1/\lambda^*\}$, where $^*$ denotes the complex conjugate. This implies that, for stable systems, all of the eigenvalues must lie in the unit circle (see for example \cite{dragt15lmf} for more details on symplectic matrices and their properties).

If we solve numerically eq. (\ref{firstordersystem}) using standard methods like say, a Runge-Kutta scheme of order $p$, to obtain $\widetilde\Phi(T)$  as an approximation to  $\Phi(T)$ we find that
\begin{equation}
	\widetilde \Phi^T\,J\,\widetilde\Phi=J + {\cal O}(h^{p})
\end{equation}
where $h$ is the time step used in the numerical integration,
and $\det \widetilde\Phi^T=1 + {\cal O}(h^{p})$ and consequently, volume preservation is not guaranteed. 
The eigenvalues of $\widetilde\Phi(T)$ will not, in general, occur in pairs and, what is even worse, in general we will find that
\[
  |\widetilde \lambda_k|=1+ {\cal O}(h^{p}) , \qquad k=1,2,\ldots,2r.
\]
Even if the exact solution is stable, the numerical solution can provide a fundamental matrix solution with eigenvalues inside and/or outside the unit circle.
Then, one should use a very small time step to avoid these undesirable effects which is similar to the step size restrictions that occur when using explicit RK methods to solve stiff equations.

For this reason, it is of great interest to study the numerical integration of the Hill's equation using symplectic integrators.


\subsection{Symplectic methods}

In the following, we consider different classes of symplectic integrators from the literature and we analyze their performance when applied to the Hill's equation. This analysis will be used to build new symplectic integrators tailored to the Hill's equation in the following sections.

\subsubsection{Implicit symplectic Runge--Kutta methods}

It is well known that standard explicit Runge--Kutta (RK) methods are not symplectic. However, the $s$-stage implicit Runge--Kutta--Gauss--Legendre (RKGL) methods are symplectic and of maximal order $2s$. RK methods are characterized by the real numbers $a_{ij}$,
$b_i$  ($i,j=1, \ldots, s$) and $c_i = \sum_{j=1}^s a_{ij}$, and for this linear problem they take the form
\begin{equation}
  \label{RK-lineal1}
\begin{split}
 Z_{i} & =  z_{n} +
 h \sum_{j=1}^{s}a_{ij}A_jZ_{j} ,
  \qquad i=1,\ldots,s    
  \\
 z_{n+1} & =  z_{n} +
 h \sum_{i=1}^{s}b_{i}A_iZ_{i},
\end{split}
\end{equation}
where $A_i=A(t_n+c_ih)$. 
This linear system can be exactly solved, however, from the computational point of view it is not advisable to use direct methods. Notice that given $Q,P\in\mathbb{C}^{r\times r}$, the computational cost to multiply these matrices is 
\[
 {\cal C}=\cost(QP)=2r^3-r^2 \ \  \mbox{flops}.
\]
We will take the cost of a method based in units of ${\cal C}$. For example, the cost to solve the system $QX=P$ is $\cost(Q^{-1}P)\simeq \tfrac43{\cal C}$. Then, na\"ive counting of multiplications results in a cost of $(2\times s)^3{\cal C}$ to solve the system (\ref{RK-lineal1}) using a direct method, which would render the method uncompetitive versus explicit methods. Obviously, this can be improved using iterative methods to solve the implicit equations (\ref{RK-lineal1}), where the cost of each iteration is only $2s{\cal C}$ (this number corresponds to the products $A_iZ_i, \ i=1,\ldots,s$ where each product, due to the sparse structure of the matrix $A$, involves two products of matrices of dimension $r\times r$). 

High order methods are useful to get accurate results while using relatively large time steps. On the other hand, large time steps can reduce that rate of convergence of the implicit methods, and in turn require more iterations. In order to preserve the symplecticity, the algorithm should be used with a very small tolerance.
One can use a fixed point iteration that - for second order equations - increases the convergence by a factor $h^2$ at each iteration. In order to preserve symplecticity to nearly round-off accuracy,  typically 5--7 iteration will be necessary.

\subsubsection{Splitting methods}

To overcome the difficulties encountered with implicit methods, symplectic splitting methods can be used in their stead. For this purpose, we express the matrix equation (\ref{Fundamental}) as an equivalent Hamiltonian system with time-dependent Hamiltonian function
\begin{equation}\label{HillsHam}
	  H(q,p,t) = \frac12 p^Tp + \frac12 q^T M(t) q
\end{equation}
with $q,p\in\mathbb{R}^{r}$. Highly efficient symplectic Runge--Kutta-Nystr\"om methods for this Hamiltonian can be found in \cite{blanes02psp}. One step from $t_n$ to $t_n+h$, of an $s$-stage method applied to solve \eqref{HillsHam} is given by
\begin{equation*}
  \begin{array}{l}
 \tau_0 := t_n,\quad  Q_0 := q_n,\quad
 P_0 := p_n,\\
  {\bf do} \ \ k=1,\ldots, s \\
     \quad Q_k  := Q_{k-1} + a_k\, h \, P_{k-1} \\
     \quad \tau_k  := \tau_{k-1} + a_k\,  h \\
     \quad P_k  := P_{k-1} - b_k\,  h\,  M(\tau_k)Q_k \\
  {\bf enddo} \\
     t_{n+1}  := \tau_s=t_n+h, \quad
     \quad q_{n+1}  := Q_s, \quad
     \quad p_{n+1} := P_s,
\end{array}
\end{equation*}
where $a_k,b_k$ are appropriate coefficients. Notice that the variable $t$ is advanced with the coefficients $a_i$ associated to the kinetic part.
This algorithm can be used for solving (\ref{Fundamental}) as follows
\begin{equation*}
	 \Phi_{n+1} =
	\begin{pmatrix}
  I & 0 \\ 
- b_s\,  h  M_s & I
\end{pmatrix}
	\begin{pmatrix}
  I & a_s hI \\ 
  0 & I
\end{pmatrix}
 \quad \cdots \quad
	\begin{pmatrix}
  I & 0 \\ 
- b_1  h  M_1 & I
\end{pmatrix}
\begin{pmatrix}  I & a_1 h I \\ 
  0 & I
\end{pmatrix}
\Phi_n,
\end{equation*}
where $\Phi_n\simeq\Phi (t_n)$ and $M_k=M(\tau_k)$. Taking into account that, in general, $\Phi_n$ is a $2r\times 2r$ dense matrix and $M_k$ is a $r\times r$ dense matrix, one step requires $s$ evaluations of the time-dependent matrix $M(t)$ (usually of low computational cost for most Hill's equations of practical interest) and $2s$ products of $r\times r$ dense matrices, i.e. the cost for one step would be $2s{\cal C}$.

In the present setting, the problem is explicitly time-dependent and highly oscillatory for most values of the parameters and symplectic methods based on the Magnus expansion can be more efficient in many cases.

\subsubsection{Magnus integrators}

The Magnus expansion \cite{magnus54ote}  expresses the solution to (\ref{Fundamental}) in the form of a single exponential
\begin{equation}\label{Magnus}
\Phi(t,0) = \exp \left( \Omega (t,0) \right) \, , \qquad 
\Omega (t,0) = \sum_{k=1}^{\infty} \Omega_{k} (t,0) ,
\end{equation}
where the first terms of the Magnus series $\{\Omega_k\}$ are given by
\begin{displaymath}
\Omega_{1} (t,0) = \displaystyle \int_{{0}}^{t} A(t_{1}) \, dt_{1} \, , \ \
\Omega_{2} (t,0) = \frac{1}{2} \displaystyle \int_{{0}}^{t}  \, \int_{t_{0}}^{t_{1}}    \left[ A (t_{1}), A (t_{2}) \right] \,        dt_{2} \, dt_{1}   , \ \ldots
\end{displaymath}
where $ \, [P,Q] = PQ - QP \, $ is the matrix commutator of $P$ and $Q$. Here $\Omega$ as well as any truncation of the series belong to the Lie algebra, and the symplectic property is preserved.

In order to obtain an approximation to $\Omega$ defined by (\ref{Magnus}) for a time step from $t_n$ to $t_{n+1}=t_n+h$, it is convenient to express it in the graded free algebra generated by $\{ \alpha_{1} , \alpha_{2}  , \ldots \}$ where
\begin{equation}\label{alpha_i}
  \alpha_{i+1} = \frac{h^{i+1}}{i!} \, \left. \frac{d^{i} A(t)}{dt^{i}} \right|_{t=t_{n}+ \frac{h}{2}},
\end{equation}
therefore $ \alpha_{i} = \mathcal{O} (h^{i})$.
The Magnus expansion $\Omega$ can be approximated to arbitrary order in this algebra.
We first consider sixth-order methods, and up to this order it suffices to take into account the Lie algebra generated by $\{ \alpha_{1}, \alpha_{2}, \alpha_{3}\}$ (see \cite{alvermann11hoc} for a relatively simple proof to this result). By Taylor expanding $A(t)$ around the midpoint $t_h+\tfrac{h}{2}$ (see \cite{blanes09tme} for a comprehensive review) it follows that we obtain a sixth-order approximation to $\Omega$ by
\begin{equation}  \label{sixth}
 \Omega^{[6]} =\alpha_{1} + \frac{1}{12} \alpha_{3} - \frac{1}{12} [12] +
  \frac{1}{240} [23] +
  \frac{1}{360} [113] - \frac{1}{240} [212] + \frac{1}{720} [1112],
\end{equation}
where $[ij\ldots kl]$ represents the nested commutator $\lbrack
\alpha_{i},[\alpha_{j}, [\ldots ,[\alpha_{k},\alpha_{l}]\ldots
]]]$ and $\Omega^{[6]} = \Omega + \mathcal{O} (h^{7}) \, .$ 
Let $c_1,\ldots,c_m$ denote the nodes of a quadrature rule of order six or higher, it is then possible to replace each $\alpha_i, \ i=1,2,3$ by a linear combination of $A_k\equiv A(t_n+c_kh), \ k=1,\ldots,m$ such that $\Omega^{[6]}$ is still an approximation to order six. Letting $b_1,\ldots,b_m$ denote the weights of the same quadrature rule, the $\alpha_i$ can be written, for example, as
\begin{equation}   \label{alpha_iA_iGeneral}
  \alpha_1 = \frac94 A^{(0)} - 15 A^{(2)},  \quad \alpha_2 = 12  A^{(1)}, \quad
  \alpha_3 = -15  A^{(0)} + 180 A^{(2)},
\end{equation}
where
\begin{equation}   \label{A_iGeneral}
   A^{(i)} = h\sum_{j=1}^m b_j\left(c_j-\frac12\right)^i A_j, 
	\qquad i=0,1,2.
\end{equation}
If we consider, for example, Gauss--Legendre (GL) collocation points where $c_i,b_i, \ i=1,2,3$ are given by:
\[
c_1=\frac{1}{2} - \frac{\sqrt{15}}{10}, \ c_2=\frac{1}{2}, \ c_3=\frac{1}{2} + \frac{\sqrt{15}}{10}, \
b_1=\frac{5}{18}, \ b_2=\frac{4}{9}, \ b_3=\frac{5}{18}
\]
and substituting into (\ref{A_iGeneral}) and (\ref{alpha_iA_iGeneral}) we have (see also \cite{blanes09tme} and references therein)
\begin{equation}   \label{alpha_iA_i}
  \alpha_1 = h A_2,  \quad \alpha_2 = \frac{\sqrt{15}h}{3} (A_3 - A_1), \quad
  \alpha_3 = \frac{10h}{3} (A_3 - 2 A_2 + A_1),
\end{equation}
where it is easy to see that $\alpha_1 =\mathcal{O} (h),  \ \alpha_2 = \mathcal{O} (h^{2}), \   \alpha_3 =\mathcal{O} (h^{3})$. Inserting this result into (\ref{sixth}) yields a straight-forward method $\exp( \Omega^{[6]} )$.

\paragraph{The computational cost of exponential of matrices}

The matrix $A(t)$ is a relatively sparse matrix and using commutators in (\ref{sixth}) reduces this sparsity. As a result, the computational cost to compute $\exp (\Omega^{[6]})$ is considerably higher than to compute $\exp (\widetilde A)$, where $\widetilde A=h\sum_i \beta_i A_i$ (with constants $\beta_i$) denotes a linear combination of $A(t)$ evaluated at different instants. This becomes obvious when examining the structure of the symplectic matrices
\begin{displaymath}
 E_1 = \exp
\begin{pmatrix}
 D & B \\
C & -D^T
\end{pmatrix},
\end{displaymath}
with $B^T=B,C^T=C$, 
versus the also symplectic
\begin{equation}\label{cheap}
  E_2 = \exp
\begin{pmatrix}
0 & I \\
C & 0
\end{pmatrix}
 \, , \qquad \text{or} \
\qquad   E_3 = \exp
 \begin{pmatrix}
0 & 0 \\
C & 0
\end{pmatrix} =
 \begin{pmatrix}
I & 0 \\
C & I
\end{pmatrix}.
\end{equation}
The exponential $E_3$ is trivial and the simple structure of $E_2$ allows to write the exponential in a simple closed form where computations can be reused and an approximation to order $2m$ can be reached at the cost of only $(m-1)\CC$ and a symplecticity can be preserved too at an extra cost of $(2+\frac13)\CC$ totaling $(m+\frac43)\CC$ (see the appendix), while this is not possible in the computation of $E_1$. 
In consequence, we find that
\[
 \cost(E_3)<\cost(E_2)\ll   \cost(E_1)
\]
and we will look for composition methods that require mostly (symplectic) exponentials of matrices with structure $E_3$ while keeping the number of exponentials of the class $E_2$ to a minimum. 
This task requires a profound analysis of the Lie algebra associated to the Hill's equation.
%

\paragraph{Commutator-free Magnus integrators}

A standard way to avoid the computation of commutators in (\ref{sixth}) which give rise to dense matrices $E_1$ is given by commutator-free (CF) Magnus integrators \cite{alvermann11hoc,free}. A sixth-order example is the following composition 
\begin{equation}\label{CF}
\exp \left( \Omega^{[6]} \right) \, 
= \, \prod_{i=1}^{s} \exp \left( \sum_{k=1}^{3} x_{i,k} \, \alpha_{k} \right) \,
= \, \prod_{i=1}^{s} \exp \left( h \gamma_i \begin{pmatrix} 0 & \delta_i I \\ C_{i} & 0 \end{pmatrix}  \right),
\end{equation}
where the coefficients $x_{i,k}$ satisfy a set of polynomial equations. Here, $C_{k} $ are linear combinations of $M(t)$ evaluated at quadrature nodes and, if $x_{i,1}\neq 0$ then $\delta_i=1, \ \gamma_i=x_{i,1}$ otherwise $\delta_i=0, \gamma_i=1$. Sixth-order methods need to take $ s \geq 5$ (a four-exponential sixth-order method exists, but it shows a very poor performance and it is not recommended in practice). In addition, some coefficients $\gamma_{i}$ are negative. 


\section{Exponential symplectic methods for the Hill's equation}

\subsection{Sixth-order methods}

The additional structure of the Hill's equation makes it possible to build new methods that improve the performance of the existing ones. 
The idea is similar to the schemes proposed for $N$th-order time-dependent linear systems \cite{bbcp}, but tailored for the Hill's equation.

The key point is to exploit the algebraic structure of $\alpha_{1}, \alpha_{2}$ and $\alpha_{3}$\footnote{Recall that the $\alpha_{i+1}$ correspond to the $i$th derivative of $A(t)$, see \eqref{alpha_i}.}, i.e., 
\begin{displaymath}
\alpha_{1} = h\begin{pmatrix}
0 & I \\ P & 0
\end{pmatrix} \,  , \ \
\alpha_{2} = h^2\begin{pmatrix}
0 & 0 \\ Q & 0
\end{pmatrix}  , \ \
\alpha_{3} = h^3 \begin{pmatrix}
0 & 0 \\ R & 0
\end{pmatrix} \,  ,
\end{displaymath}
where, in the case of the GL quadrature rule (\ref{alpha_iA_i}) we have that
\[
   P = -M_2,  \quad h Q = -\frac{\sqrt{15}}{3} (M_3 - M_1), \quad
  h^2R = -\frac{10}{3} (M_3 - 2 M_2 + M_1).
\]
Since $\alpha_{2}$ and $\alpha_{3} $ are nilpotent matrices of degree two, they can be exponentiated trivially,
\[
  \exp(x\alpha_{2}+y\alpha_{3}) = I + x\alpha_{2}+y\alpha_{3}, \quad x,y\in\mathbb{C}
\]
and the exponential of $\alpha_{1}$ has a considerably lower computational cost than for a full matrix. In addition, we observe that 
\[
 [\alpha_{2},\alpha_{3}]=[23]=0 
\qquad \mbox{and} \qquad
 [212] =  h^5\begin{pmatrix}
0 & 0 \\ 2Q^2 & 0
\end{pmatrix},
\]
i.e., $[212]$ is also nilpotent with similar structure to $\alpha_{2}$ and $\alpha_{3}$.
Since the exponentials of matrices that contain $\alpha_1$ will be the most costly part of the new schemes, we analyze new compositions by the number of such exponentials involved.


\paragraph{One-$\alpha_1$-exponential method}
In \cite{bbcp}, the following fourth-order composition was proposed
\begin{equation}\label{CFM41}
   \Phi_1^{[4]} = \exp\left( \frac1{24} (2\alpha_2 + \alpha_3) \right) \,
  \exp\left(\alpha_1 \right)  \, \exp\left( \frac1{24} (-2\alpha_2 + \alpha_3)\right).
\end{equation}
As we have just established, $[212]$ can be added to the first and last exponentials without significantly increasing the cost and furthermore
\[
 [1112]= \begin{pmatrix}
W & 0 \\ 0 & -W^T
\end{pmatrix},
\]
with $W=h^5(3QP+PQ)$ is very small and its exponential can be approximated with two products and one inversion for the symplectic case ($QP=(PQ)^T$), i.e., cost $(2+\tfrac43){\cal C}$ by
\[
 \exp \big([1112]\big)=
\begin{pmatrix}  \Lambda & 0 \\ 0 & \Lambda^{-T} 
	\end{pmatrix}
+ {\cal O}(h^{15}).
\]
where $\Lambda =I+W+\frac12W^2$. The computation of the inverse ensures symplecticity. 
Using these observations, we propose the composition
\begin{multline}\label{Magnus 61}
   \Phi_1^{[6]} = 
	\exp\big( x_6[1112] \big) \,
	\exp\big( x_3\alpha_2 + x_4\alpha_3 + x_5[212]\big) \, 
	  \exp\big(x_1 \alpha_1 + x_2 \alpha_3 \big)\\
\times \, 	\exp\big(-x_3\alpha_2 + x_4\alpha_3 + x_5[212]\big) 
	\exp\big( x_6[1112] \big).
\end{multline}
There is only one solution for the coefficients $x_i$ given by
\[
x_1=1,\quad x_2=\frac{1}{20},\quad 
x_3=\frac{1}{12},\quad  x_4=\frac{1}{60},\quad 
x_5=-\frac{1}{2880},\quad  x_6=\frac{1}{720}.
\]
It is important to remark that, while the cost of $\exp\big( x_6[1112] \big)$ is $(2+\tfrac43){\cal C}$, in the flow of the algorithm this matrix is multiplied with a dense matrix $\Phi$ (i.e. $\exp\big( x_6[1112] \big)\Phi$) which increases the total cost to $(2+\tfrac43){\cal C}+4{\cal C}$. The computational effort can be reduced by generalizing the (first same as last) FSAL property: We concatenate two steps, the last exponential of one step and the first one in the following one\footnote{For the usual FSAL property, the two exponentials are identical and can be concatenated without error.}, and since the exponents are small, we commit an error $e^{h^5V}e^{h^5W}=e^{h^5(V+W)+{\cal O}(h^{10})}$ and count the cost of only one exponential  $\exp\big( x_6[1112] \big)$ per step.

The total cost per step is $27+\frac23$, 
$$\cost(\Phi_1^{[6]})
=
\underbrace{6+\frac43}_{\cost(\exp([1112])\Phi}+\underbrace{\vphantom{\frac{1}{2}}(1+2)+2}_{2\,\cost(\exp(\alpha_2+\alpha_3+[212]))\Phi}+\underbrace{\vphantom{\frac{1}{2}}(7+\frac13)+8}_{\cost(\exp(\alpha_1+\alpha_2+\alpha_3))\Phi},
$$
where the $\cost(\exp(\alpha_2+\alpha_3+[212]))\Phi=1+2$ comes from one product for [212] (is computed only once and stored to be used in the second one) and a matrix multiplication with a sparse matrix and $\cost(\exp(\alpha_1+\alpha_2+\alpha_3))\Phi=(7+\frac13)+8$ stems from an approximation to order 13 in \eqref{Exp} (5 products) and preservation of symplecticity  (2+$\frac13$ products).

\paragraph{Two-$\alpha_1$-exponential method}
The next composition has enough parameters to reach order six,
\begin{multline}\label{Magnus 64}
   \Phi_2^{[6]} = 
\exp \left(x_{1} \alpha_{2} + x_{2} \alpha_{3} + x_{3} \left[\alpha_{2},\alpha_{1},\alpha_{2} \right] \right) \, 
\exp \left(
x_{4} \alpha_{1} + x_{5} \alpha_{2} + x_{6} \alpha_{3}
\right) \\
\times \, \exp \left(
x_{4} \alpha_{1} - x_{5} \alpha_{2} + x_{6} \alpha_{3}
\right) \, 
\exp \left(
-x_{1} \alpha_{2} + x_{2} \alpha_{3} + x_{3} \left[ \alpha_{2},\alpha_{1}, \alpha_{2} \right]
\right) .
\end{multline}
There is only one solution for the coefficients $x_i$,
\[
x_1=\frac{1}{60},\quad x_2=\frac{1}{60},\
x_3=\frac{1}{43200},\hspace*{1.3mm} x_4=\frac{1}{2},\
x_5=\frac{2}{15},\ \hspace*{5.8mm}x_6=\frac{1}{60},
\]
and the composition takes the form
\begin{displaymath}
 \Phi_2^{[6]} = \begin{pmatrix} I & 0 \\ hC_{4} & I \end{pmatrix}
\exp \left(\frac{h}2 \begin{pmatrix} 0 & I \\ C_{3} & 0 \end{pmatrix}\right)
\exp \left(\frac{h}2\begin{pmatrix} 0 & I \\ C_{2} & 0 \end{pmatrix}\right)
\begin{pmatrix} I & 0 \\ hC_{1} & I \end{pmatrix},
\end{displaymath}
where $ C_{i}$, $ \, i=1,2,3,4 \, $ are linear combinations of $M(t) $ evaluated in a set of quadrature points and $C_1, C_4$ additionally contain one product of such linear combinations. If the sixth-order Gaussian quadrature rule is used, the method given in Table~\ref{method62} is obtained, but any other quadrature rule of order six or higher can also be used.

The total cost preserving symplecticity and up to order 13 ($7+\frac13$ products) in \eqref{Exp} is
$$\cost(\Phi_2^{[6]})
=
\underbrace{\vphantom{2((7+\frac13)+8)}(1+2)}_{\cost(\exp(\alpha_2+\alpha_3+[212])))\Phi}+\underbrace{2((7+\frac13)+8)}_{2\,\cost(\exp(\alpha_1+\alpha_2+\alpha_3))\Phi}=33+\frac23.
$$
Here, we exploited that the last exponential can be concatenated with the first one in the following step at no extra cost, and it is not counted, a property which we call \textit{first commutes with last} (FCWL).

\paragraph{Three-$\alpha_1$-exponential method}
We analyze now the following composition with three exponentials
\begin{multline}\label{Magnus63}
   \Phi_3^{[6]} = 
\exp \left( x_{1} \alpha_{1} + x_{2} \alpha_{2} + x_{3} \alpha_{3} \right) \,
\exp \left( x_{4} \alpha_{1} + x_{5} \alpha_{3} + x_{6} \left[ \alpha_{2} , \alpha_{1} , \alpha_{2} \right] \right) \, \\
\times \, \exp \left( x_{1} \alpha_{1} - x_{2} \alpha_{2} + x_{3} \alpha_{3} \right) \, .
\end{multline}
There are only two solutions, one of them with  $x_1, x_4>0$ given by:
	\begin{alignat}{3}
x_1&=\frac{1}{10}(5-\sqrt{5}),\qquad x_2=\frac{1}{24}(5-\sqrt{5}),\qquad 
x_3=\frac{1}{60}(5-\sqrt{5}),\\ 
x_4&=\frac{1}{\sqrt{5}},\qquad
x_5=\frac{1}{60}(-5+2\sqrt{5}),\ \hspace*{1.3mm}x_6=\frac{1}{8640}(-11+5\sqrt{5}) .
	\end{alignat}
This composition has the structure
\begin{displaymath}
\exp \left( \gamma_{3} h \begin{pmatrix} 0 & I \\ C_{3} & 0 \end{pmatrix}\right) \, 
\exp \left( \gamma_{2} h \begin{pmatrix} 0 & I \\ C_{2} & 0 \end{pmatrix}\right) \,
\exp \left( \gamma_{1} h  \begin{pmatrix}0 & I \\ C_{1} & 0 \end{pmatrix} \right) \, ,
\end{displaymath}
where $ \gamma_{1}= \gamma_{3}=x_1, \  \gamma_{2}=x_4$.

The total cost is
$$\cost(\Phi_3^{[6]})
=
\underbrace{\vphantom{2(4+8)}1+(7+\frac13)+8}_{\cost(\exp(\alpha_1+\alpha_3+[212])))\Phi}+\underbrace{2((7+\frac13)+8)}_{2\,\cost(\exp(\alpha_1+\alpha_2+\alpha_3))\Phi}=47.
$$

\subsection{Eigth-order methods}

We extend the previous analysis to build eight-order methods. 
The Magnus expansion up to order eight for a general linear system using the Lie algebra generated by $\{\alpha_{1},\alpha_{2}, \alpha_{1},\alpha_{4} \}$ reads
\begin{eqnarray*}  \label{eighth}
 \Omega^{[8]} & = & \alpha_{1} + \frac{1}{12} \alpha_{3} - \frac{1}{12} [12] 
  + \frac{1}{240} [23] +
  \frac{1}{360} [113] - \frac{1}{240} [212] + \frac{1}{720} [1112]  \\
	& & - \frac{1}{80}[14]-\frac{1}{1344} [34] -\frac{1}{2240}[124]
	+\frac{1}{6720}[223] +\frac{1}{6048}[313]-\frac{1}{840}[412]   \\
&&
+\frac{1}{6720}[1114]-\frac{1}{7560}[1123] +\frac{1}{4032}[1312]+\frac{11}{60480}[2113]-\frac{1}{6720}[2212]   \\
&&
-\frac{1}{15120}[11113]-\frac{1}{30240}[11212]  +\frac{1}{7560}[21112] -\frac{1}{30240}[111112]. \nonumber
\end{eqnarray*}
Further significant simplifications occur for the Hill's problem since the following commutators cancel
\[
 [23]=[34]=[24]=[223]=[1123]=[2212]=0.
\]
In addition, the matrices $ [212], \ [313], \ [412] $
are nilpotent and these elements can be included in the exponents at a minor extra cost. We can also use other nilpotent matrices of even order, such as $[312]$, if they are distributed into the composition skew-symmetrically (the element $[213]$ is not included because of the Jacobi identity $[123]+[231]+[312]=0$ and $[123]=0$ so that 
$[312]=[213]$).


We consider a number of compositions with enough independent parameters to solve the order conditions using four and five exponentials. In each case, there are many different possible compositions leading in some cases to complex-valued solutions, to one or several real solutions or to families of solutions in terms of free parameters. In the following we present the method that provided the best performance in practice among the methods studied.

\paragraph{Five-$\alpha_1$-exponential method}
The following composition
\begin{multline*}\label{}
\Phi_5^{[8]}=
\exp   \big(x_7 \alpha_1+x_8\alpha_2+x_9 \alpha_3+x_{10} \alpha_4+x_{11} [212]\big)  \\
\exp   \big(x_{12} \alpha_2+x_{13} \alpha_3+x_{14} \alpha_4+x_{15} [212]
	+x_{16} [313]\big)  \\
\exp   \big(x_3 \alpha_1+x_4 \alpha_2+x_5 \alpha_3+x_6 \alpha_4\big)  \\
\exp   \big(x_{1} \alpha_1+ x_{2}\alpha_3\big)  \\
 \exp  \big(x_3 \alpha_1-x_4 \alpha_2+x_5 \alpha_3-x_6 \alpha_4\big)  \\
 \exp  \big(-x_{12} \alpha_2+x_{13} \alpha_3-x_{14} \alpha_4+x_{15} [212]
	+x_{16} [313]\big)\\
	\exp   \big(x_7 \alpha_1-x_8 \alpha_2+x_9 \alpha_3-x_{10} \alpha_4+x_{11} [212]\big)
\end{multline*}
has two real valued solutions at order six, and the one with smallest coefficients is
{\footnotesize\[
\begin{array}{lllll}
x_1    &= 0.6403363286379515,			& x_2 &= 0.0433501199827269,\\
x_3    &=-0.4017895263297271, 		& x_4 &=-0.1170180583697493,\\
x_5    &=-0.1038563759039891, 		&x_6  &=-0.0376728349617945\\
x_7    &= 0.5816213620107513,     &	x_8 &= 0.2609350592183406, \\
x_9    &= 0.1157777422250884	  ,&x_{10}&= 0.0506748377294480,\\ 
x_{11} &=-0.0000936846387697,		 &x_{12}&=-0.0127292796833454,\\
x_{13} &= 0.0080702403542039, 	 &x_{14}&=-0.0017487133111753,\\
x_{15} &=-0.0000928250351798,		 &x_{16}&= 0.0001835812673590  
\end{array}
\]	}
The computational cost of this method is: 2${\cal C}$ to compute $[212],[313]$, $2\times 2$${\cal C}$ to compute the products of nilpotent matrices, and $5\times ((7+\frac13)+8)$${\cal C}$ to approximate the exponentials up to 12th-order, so the total cost is $82+\frac23 {\cal C}$.

The computational cost of the different 8th-order methods we considered is about 70-85 ${\cal C}$ where the exponentials are approximated up to 12th order), it is approximately half the cost of the commutator-free methods of the same order obtained in \cite{alvermann11hoc} which require 11 exponentials, that is $168+\frac23 {\cal C}$.

For the implementation of the eighth-order methods, the derivatives $\alpha_i$ are replaced by momentum integrals using the substitution rules
\begin{equation}  \label{eq.3.9}
 \begin{array}{ccl} \displaystyle
  \alpha_1 & = & \frac{3}{4} \, ( 3 A^{(0)} - 20 A^{(2)}),   \\
  \alpha_2 & = & 15 \, (5 A^{(1)} - 28 A^{(3)}),  
 \end{array}
\qquad
 \begin{array}{ccl}  \displaystyle
  \alpha_3 & = & -15 \,( A^{(0)} - 12 A^{(2)}),   \\
  \alpha_4 & = & -140 \,(3 A^{(1)} - 20 A^{(3)}),
 \end{array}
\end{equation}
and the $A^{(i)}, \ i=1,2,3,4$ are approximated with a standard rule following (\ref{A_iGeneral}).
For instance, using the 8th order Gauss-Legendre quadrature rule, we obtain
\[
c_i=\frac12-v_i, \ 
c_{5-i}=\frac12+v_i, \quad
b_i=\frac12-w_i, \ 
b_{5-i}=\frac12+w_i, 
\]
$i=1,2$, and where
\[
  v_1 = \frac12  \sqrt{\frac{3+2 \sqrt{6/5}}{7}},   \ \
  v_2 = \frac12  \sqrt{\frac{3-2 \sqrt{6/5}}{7}},\quad
   w_1 = \frac{1}{2} - \frac{1}{6} \sqrt{\frac{5}{6}}, \ \ 
   w_2 = \frac{1}{2} + \frac{1}{6} \sqrt{\frac{5}{6}}.
\]
Letting $S_1=A_1+A_4$, $S_2=A_2 + A_3$, $R_1=A_4-A_1$, $R_2=A_3-A_2$, we reach the expression
\begin{equation}  \label{}
\begin{array} {ll}
	\displaystyle A^{(0)} = \frac{1}{2} \left(w_1 S_1 +  w_2 S_2\right), & 
	\displaystyle	A^{(1)} = \frac{1}{2} \left(w_1v_1 R_1 +  w_2v_2  R_2\right), \\[4mm]
	\displaystyle	A^{(2)} = \frac{1}{2} \left(w_1v_1^2 S_1 +  w_2v_2^2  S_2\right), \qquad & 
	\displaystyle	A^{(3)} = \frac{1}{2} \left(w_1v_1^3 R_1 +  w_2v_2^3  R_2\right).	
\end{array}
\end{equation}
%

\section{Numerical examples}
We now study the performance of the new methods on the numerical integration of different Hill's equations. The number of $r\times r$ matrix-matrix products, $k$, to compute the product $\Psi \Phi$, where $\Psi$ denotes the method and $\Phi$ the fundamental matrix solution, is given in parenthesis as $k\, \CC$, and it is taken as the cost of the method. We use 12th-order symplectic approximations to evaluate the exponentials $E_2$ at the cost of  $(7+\frac13)\, \CC$ per exponential. The following methods are considered:
\begin{itemize}
	\item RK$_7^{[6]}$: A 7-stage explicit (non symplectic) RK method that only requires 3 new evaluations of $A(t)$ per step given in \cite[p. 203-205]{butcher87tna} (14  ${\cal C}$).
	\item  RKGL$^{[6]}$: The 3-stage implicit symplectic RKGL method (we count an average of 6 iterations per step for a total cost of 36  ${\cal C}$).
	\item  RKN$_{11}^{[6]}$: The 11-stage explicit symplectic RKN method given in \cite{blanes02psp} (22  ${\cal C}$). This method requires 11 new evaluations of $A(t)$ per step that are not counted into the cost. 
	\item $\Phi_5^{[6]}$: The five-exponential commutator-free Magnus integrator from \cite{free} ($(76+\frac23)  {\cal C}$). 
	\item $\Phi_i^{[6]}$: The new $ \Phi_i^{[6]}, \ i=1,2,3$ methods ($(27 + \frac23)\CC, (33+ \frac23)\CC, 47{\cal C}$, respectively).
	\item $\Phi_5^{[8]}$: The new  $ \Phi_5^{[8]}$ method ($(82+ \frac23){\cal C}$).
\end{itemize}

\subsection{The Mathieu equation}
As a test bench to compare the performance of the methods we employ Mathieu's equation
\begin{displaymath}
x'' + \left( \omega^{2} + \varepsilon \cos (2t) \right) x = 0 .
\end{displaymath}
In a first experiment to test the qualitative behavior, we compute the fundamental matrix solution for one period (at $T=\pi$) with the identity matrix as the initial conditions using the time step $h=\pi/10$ for a set of values $\omega$ (parametric resonances occur about $\omega=k, \ k\in\mathbb{N}$. 
We let $\varepsilon=5$ and $\omega=j\Delta \omega=j\frac1{200},\ j=0,1,\ldots,1020$, and consider the following methods: RK$_7^{[6]}$, RKGL$^{[6]}$, RKN$_{11}^{[6]}$ and $\Phi_2^{[6]}$. 

\begin{figure}[H]
  \begin{subfigure}[t]{\textwidth}
	\tikzsetnextfilename{fig1a}
	{\footnotesize\input{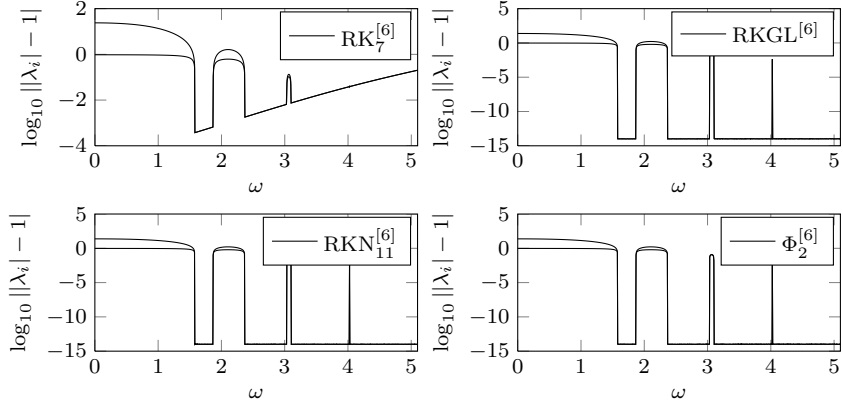}}
  \caption{Each line corresponds to one of the two eigenvalues.}
  \label{figure1a}
\end{subfigure}
\begin{subfigure}[t]{\textwidth}
	\tikzsetnextfilename{fig1b}
	{\footnotesize\input{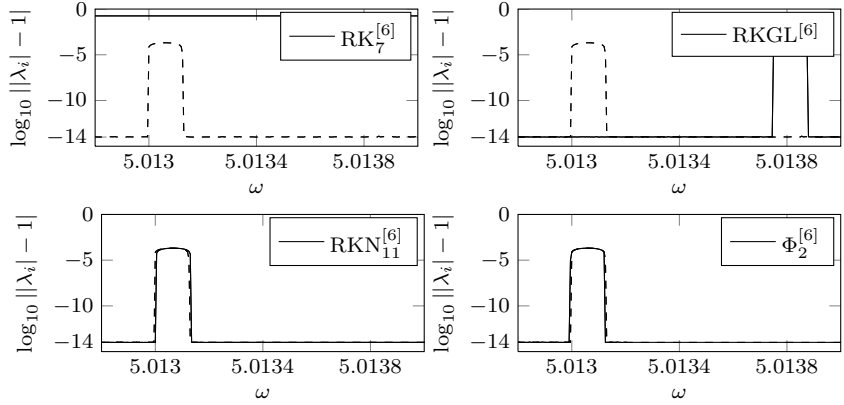}}
  \caption{Zoom of Figure~\ref{figure1a} about $\omega=5$. The dashed line correspond to the reference solution obtained numerically to sufficiently high accuracy}
  \label{figure1b}
	\end{subfigure}
	\caption{Distance of the eigenvalues to the unit circle versus $\omega$ of the explicit non-symplectic method RK$_{7}^{[6]}$, and the symplectic methods the RKGL$^{[6]}$, RKN$_{11}^{[6]}$ and $\Phi^{[6]}_2$.}
\end{figure}

We compute the eigenvalues and their distance to the unit circle,
$  |\lambda_1|-1 $ and $  |\lambda_2|-1.$
The results, in logarithmic scale\footnote{We compute $\log\big||\lambda_i|-1+\delta\big|$ with $\delta=10^{-14}$ to avoid singularities.}, are shown in Figure~\ref{figure1a} (a zoom about $\omega=5$, taking smaller values of $\Delta \omega$, is shown in Figure~\ref{figure1b} together with the reference solution, dashed lines, that was computed numerically to very high accuracy). 
From the results, it is clear that preservation of symplecticity is crucial. The non symplectic RK$_7^{[6]}$ method is about 1.5--2.5 times faster than the symplectic ones for the same time step, but requires a much smaller time step to reach similar accuracy.

\setlength\figurewidth{.35\textwidth}%
\setlength\figureheight{.15\textwidth}

\begin{figure}[htp!]
  \setlength\figurewidth{.9\textwidth}%
	\setlength\figureheight{.35\textwidth}
	\tikzsetnextfilename{fig2}
	{\pgfplotsset{every axis/.style={mark repeat=10}}
	{\footnotesize\input{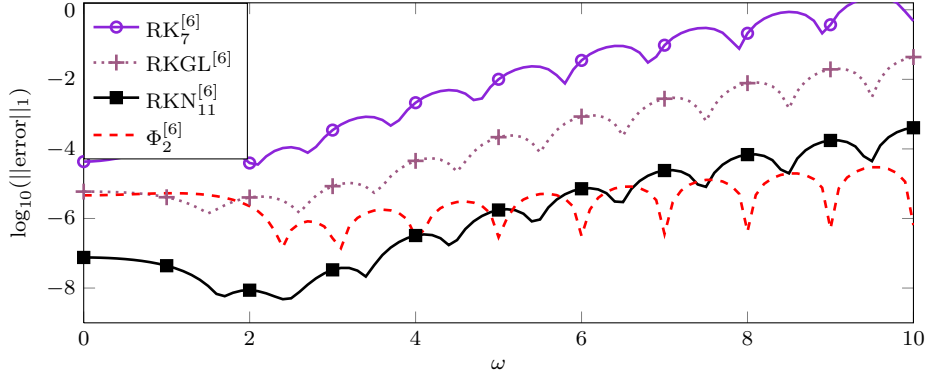}}}
  \caption{Error of the fundamental matrix solution at one period of integration versus $\omega$ with $\varepsilon=5$.}
  \label{figure2}
\end{figure}

To gain insight on how the accuracy depends on the frequency of the oscillatory solution, we repeat the same numerical experiment for $h=\pi/20$ and measure the $L_1$-norm of the error in the fundamental matrix solution (where the reference solution was obtained with a sufficiently small time step). The results are shown in Figure~\ref{figure2}. We observe that the new exponential method shows a smaller error growth with $\omega$. As a result, the same time-dependent function evaluations can be used for many different values of $\omega$ in the exponential method while the RKGL$^{[6]}$  and RKN$_{11}^{[6]}$ methods should be used with smaller time steps as $\omega$ increases and the time-dependent functions need to be recalculated.

\begin{figure}[htp!]
  \centering
		\setlength\figurewidth{.9\textwidth}%
	\setlength\figureheight{.35\textwidth}
	\tikzsetnextfilename{fig3}
	{\footnotesize\begin{tikzpicture}

\begin{axis}[%
width=\figurewidth,
height=\figureheight,
scale only axis,
xtick={2,2.2,2.4,2.6,2.8,3}, 
ytick={-10,...,0},
xmin=2,
xmax=2.95,
xlabel={ $\log_{10}(\cost)$},
ymin=-8.5,
ymax=-1,
ylabel={$\log_{10}(\|\text{error}\|_1)$},
legend style={at={(0.001,0.008)},anchor=south west,draw=black,fill=white,legend cell align=left}
]
\addplot [styleRKexplicit]
table[row sep=crcr]{
1.6232492903979	6.08593460677624	\\
1.84509804001426	2.45142748262425	\\
2.04921802267018	0.333481088445364	\\
2.22530928172586	-0.87880255960005	\\
2.40140054078154	-1.72849086257064	\\
2.57749179983723	-2.71699456781513	\\
2.75891189239797	-3.78194878654108	\\
2.93851972517649	-4.85320943839304	\\
3.11461098423217	-5.90929481545277	\\
3.29225607135648	-6.97648107075864	\\
3.46834733041216	-8.03460706420063	\\
3.64443858946784	-9.09244138836557	\\
3.82098917641605	-10.1517473048571	\\
3.99738638439731	-11.198446301704	\\
4.17347764345299	-12.1341245536316	\\
};
\addlegendentry{RK$_7^{[6]}$};

\addplot [styleRKGL]
table[row sep=crcr]{
2.03342375548695	0.658378443714222	\\
2.25527250510331	-0.291179350452254	\\
2.45939248775923	-1.40602449801055	\\
2.63548374681491	-2.42516934221367	\\
2.81157500587059	-3.46533268652733	\\
2.98766626492627	-4.51465333993003	\\
3.16908635748702	-5.59991106432078	\\
3.34869419026554	-6.67615470914626	\\
3.52478544932122	-7.73210083536741	\\
3.70243053644553	-8.79774943075767	\\
3.87852179550121	-9.85476192047971	\\
4.05461305455689	-10.9146406111463	\\
4.2311636415051	-12.0127190069438	\\
4.40756084948636	-13.0612355236735	\\
4.58365210854204	-12.7298261227014	\\
};
\addlegendentry{RKGL$^{[6]}$};

\addplot [styleRKN
]
table[row sep=crcr]{
1.81954393554187	-0.75936694352134	\\
2.04139268515822	-1.73273992854772	\\
2.24551266781415	-3.40561909263483	\\
2.42160392686983	-4.47664159631792	\\
2.59769518592551	-5.53942609251657	\\
2.77378644498119	-6.59873764851224	\\
2.95520653754194	-7.68851268613038	\\
3.13481437032046	-8.76673836635381	\\
3.31090562937614	-9.82399078606793	\\
3.48855071650044	-10.8936208754268	\\
3.66464197555613	-11.9935656258309	\\
3.84073323461181	-12.988998590429	\\
4.01728382156002	-12.6191130638544	\\
4.19368102954128	-12.6242362492335	\\
4.36977228859696	-12.9370021718783	\\
};
\addlegendentry{RKN$_{11}^{[6]}$};

\addplot [stylePhiFiveSix
]
table[row sep=crcr]{
2.36172783601759	-1.59393907097173	\\
2.58357658563395	-3.49427217795127	\\
2.78769656828987	-5.97755343551674	\\
2.96378782734556	-8.31678044135805	\\
3.13987908640124	-9.3083147813016	\\
3.31597034545692	-10.3443875441392	\\
3.49739043801767	-11.3934524081065	\\
3.67699827079618	-12.2360037529276	\\
3.85308952985187	-11.7721746719057	\\
4.03073461697617	-12.2459126502904	\\
4.20682587603185	-11.6933748697908	\\
4.38291713508753	-11.1663364184455	\\
4.55946772203574	-11.4933316519389	\\
4.73586493001701	-10.7638810737066	\\
4.91195618907269	-10.2089974265851	\\
};
\addlegendentry{$\Phi^{[6]}_5$};

\addplot [stylePhiOneSix
]
table[row sep=crcr]{
1.88649072517248	-0.872480005395293	\\
2.10833947478884	-1.28570874652442	\\
2.31245945744476	-4.28017428636766	\\
2.48855071650044	-5.08279320007582	\\
2.66464197555613	-6.26772792459273	\\
2.84073323461181	-7.37218738291783	\\
3.02215332717255	-8.47977172561632	\\
3.20176115995107	-9.56473786231236	\\
3.37785241900675	-10.6206067283094	\\
3.55549750613106	-11.6173299183042	\\
3.73158876518674	-12.4801526805295	\\
3.90768002424242	-12.376350016917	\\
4.08423061119063	-11.6362651932938	\\
4.26062781917189	-12.0106186514051	\\
4.43671907822758	-12.4187558873421	\\
};
\addlegendentry{$\Phi^{[6]}_1$};

\addplot [stylePhiTwoSix
]
table[row sep=crcr]{
2.00432137378264	-0.846700004590116	\\
2.226170123399	-3.22846969319252	\\
2.43029010605492	-5.15339199027425	\\
2.60638136511061	-6.87492429345885	\\
2.78247262416629	-7.95984743706447	\\
2.95856388322197	-9.02710258350111	\\
3.13998397578272	-10.1193584522673	\\
3.31959180856123	-11.1913189326156	\\
3.49568306761692	-12.1516782124718	\\
3.67332815474122	-12.8897633953249	\\
3.8494194137969	-11.9113460304539	\\
4.02551067285258	-12.0384509639585	\\
4.20206125980079	-11.9005701348605	\\
4.37845846778206	-11.6830563578883	\\
4.55454972683774	-11.0000464501477	\\
};
\addlegendentry{$\Phi^{[6]}_2$};

\addplot [stylePhiThreeSix]
table[row sep=crcr]{
2.14921911265538	-1.00715444103825	\\
2.37106786227174	-3.26565227396244	\\
2.57518784492766	-5.62719558059684	\\
2.75127910398334	-7.44420066615289	\\
2.92737036303902	-8.49132092687046	\\
3.1034616220947	-9.54566717341792	\\
3.28488171465545	-10.6312344922056	\\
3.46448954743397	-11.6191716774744	\\
3.64058080648965	-12.5173174520906	\\
3.81822589361396	-11.8731338798345	\\
3.99431715266964	-11.9370132501293	\\
4.17040841172532	-11.7713764418648	\\
4.34695899867353	-11.9123693835212	\\
4.52335620665479	-10.3623725946053	\\
4.69944746571047	-10.8977516515188	\\
};
\addlegendentry{$\Phi^{[6]}_3$};

\addplot [stylePhiFiveEight
]
table[row sep=crcr]{
2.39967372148104	-1.77575701456974	\\
2.62152247109739	-4.71021474537694	\\
2.82564245375332	-5.85789216070236	\\
3.001733712809	-7.29698400972339	\\
3.17782497186468	-8.69283757818043	\\
3.35391623092036	-10.0952634901302	\\
3.53533632348111	-11.5004767316635	\\
3.71494415625963	-12.1855068282324	\\
3.89103541531531	-12.8545694299905	\\
4.06868050243961	-12.7249119266255	\\
4.2447717614953	-12.1238417584583	\\
4.42086302055098	-11.5144521571229	\\
4.59741360749919	-11.7361720281774	\\
4.77381081548045	-11.4087634888775	\\
4.94990207453613	-10.3067319432012	\\
};
\addlegendentry{$\Phi^{[8]}_5$};

\end{axis}
\end{tikzpicture}%
}
  \caption{Error at the final time $t=\pi$ versus the number of products $\mathcal{C}$ in double logarithmic scale.}
  \label{fig3}
\end{figure}
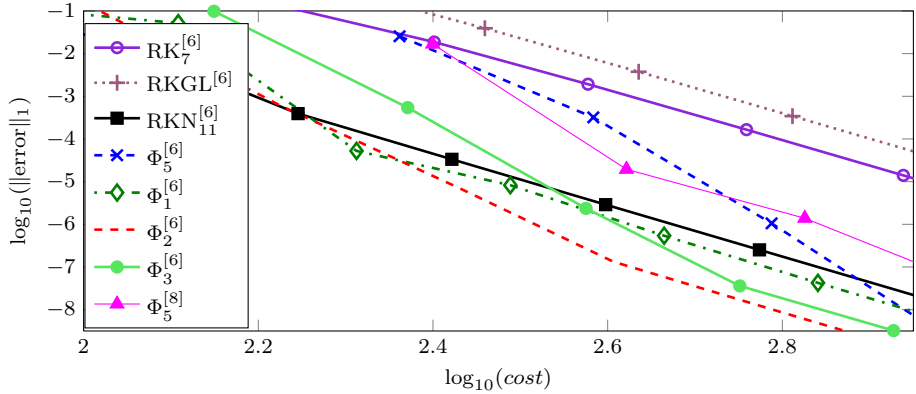

We now analyze the efficiency of the integrators. We choose a moderate value of $\omega=5 $, and $\varepsilon = 1 $ and integrate for $ \, t \in [0,\pi]$.  We measure the $L_1$-norm of the error of the fundamental matrix at the final time versus the computational cost (in units of $\mathcal{C}$). The results are shown in Figure~\ref{fig3}, where we observe that the new sixth-order method with two exponentials is superior for all accuracies.

For very small values of $\omega$ (non oscillatory problems) the symplectic RKN$_{11}^{[6]}$ integrator is the method of choice but as $\omega$ grows and the system becomes oscillatory, the new methods show the best performance.


\subsection{Matrix Hill's equation}
Let us now consider the following matrix Hill's equation
\begin{displaymath}
x'' + \left( A + B_1 \cos (2t) + B_2 \cos (4t) \right) x = 0 
\end{displaymath}
with $A,B_1,B_2\in\mathbb{R}^{r\times r}$. We take $A=r^2I+D$ where $D$ is the Pascal matrix
\[
D_{1i}=D_{i1}=1, \quad  i=1,\ldots,r, \qquad  
D_{ij}=D_{i-1,j}+D_{i,j-1}, \quad 1<i,j\leq r,
\]
which is a symmetric positive definite dense matrix. We set $B_1=\varepsilon I, \ B_2=\frac{1}{10}\varepsilon I$, $\varepsilon =r$ and $\varepsilon =\frac{1}{10}r$  and compute the solution for $r=5$ and $r=7$.

\begin{figure}[htp!]
  \centering
	\setlength\figurewidth{.38\textwidth}%
	\setlength\figureheight{.4\textwidth}
	\tikzsetnextfilename{fig4}
	{\footnotesize\input{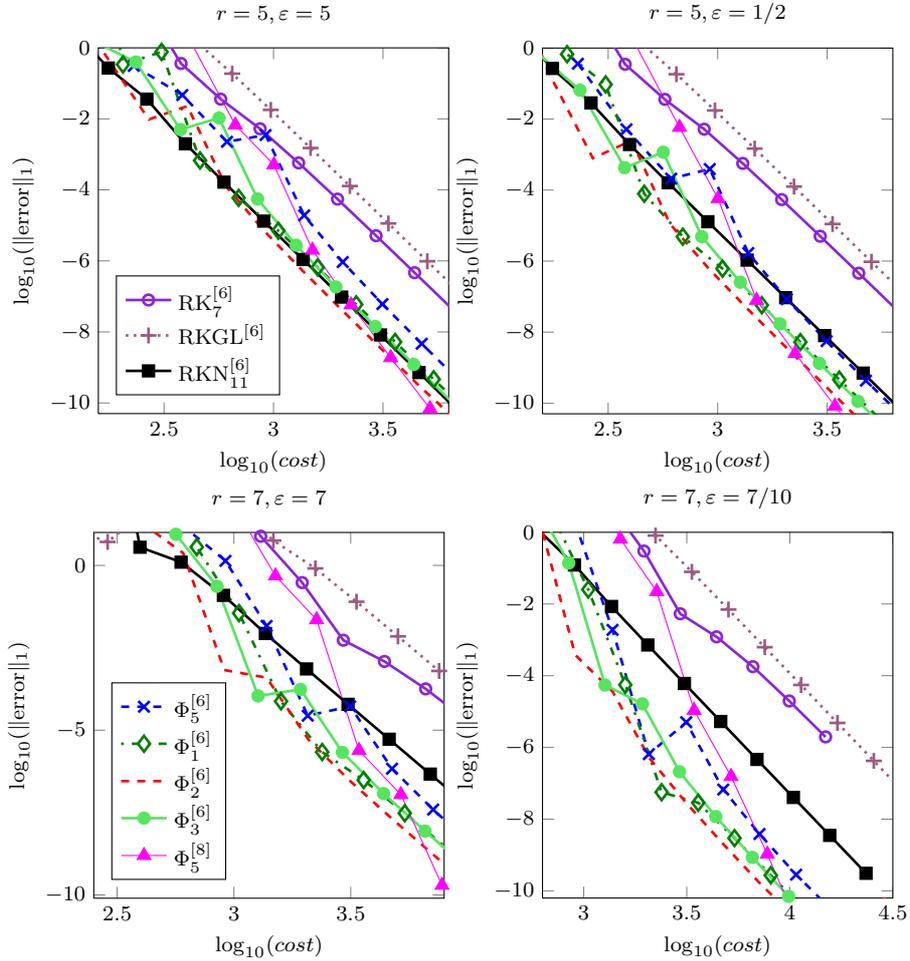}}
  \caption{Efficiency plot for Example~2. The legend has been split over two panels.}
  \label{fig4}
\end{figure}

The error of the fundamental matrix solution is measured in the $L_1$-norm. The results are shown in Figure~\ref{fig4}. 

We clearly observe that as $r$ increases (oscillatory problem) as well as $\varepsilon$ decreases (approaching the autonomous problem) the new exponential methods show the best performance. The eighth-order method is only superior when very accurate results are desired being an open problem to know if there exist other composition leading to more efficient methods that are superior for medium to high accuracy.

\section{Conclusions}

We have studied the numerical integration of the matrix Hill's equation using methods that accurately reproduce the parametric resonances of the exact solution. We are mainly interested in the Hamiltonian case which is the most frequent one in practice, namely when the Hill's equations originate from a Hamiltonian function. In this case the fundamental matrix solution is  a symplectic matrix and we illustrate the importance of the preservation of this property by the numerical integrators. 

We have presented new symplectic sixth- and eighth-order symplectic exponential integrators that are tailored to the matrix Hill's equation. Exponential integrators usually show very accuracy for stiff and oscillatory problems but at a relatively high computational cost that made them uncompetitive versus existing explicit sympletic Runge-Kutta-Nystr\"om methods when applied to the matrix Hill's equation. However, we show that a class of matrix exponentials can be very efficiently approximated while preserving the symplectic structure, and we have built new families of methods based on exponentials of this type. Several sixth- and eighth-order methods using compositions of one to five exponentials are considered. The numerical experiments showed the high performance of the new methods. Among the methods obtained, a sixth-order two-exponential method showed the best performance. The eighth-order methods obtained using different compositions had large coefficients that turned into relatively large truncation errors and, in addition, required a higher order truncation for the approximation to the exponentials. It is left as an open question if different eighth-order compositions with small coefficients exist which could show a higher performance.

The methods obtained in this work can also be used for the numerical integration of the weekly damped Hill's equation with a week non-linear interaction as the basic method to solve perturbed time-dependent linear problems as shown in \cite{bader11fmf,blanes10sac,seydaoglu14hos}. 

\appendix

\section{Efficient symplectic approximation of $E_2$}

We seek a symplectic approximation to 
\begin{equation}\label{Exp}
E_2=\exp \left( \gamma h 
\begin{pmatrix} 0 & I \\ C & 0 \end{pmatrix} \right) =\exp \left( \tau\begin{pmatrix}0 & I \\ C & 0 \end{pmatrix} \right)=\begin{pmatrix} 
\sigma & 
\mu \\ 
-C\mu & 
\sigma   \end{pmatrix} \, 
\end{equation}
with $\tau=\gamma h$ 
%
and
\begin{displaymath}
\sigma =\sum_{n=0}^{\infty} \frac{\tau^{2n}}{(2n)!} C^n, \qquad
\mu =\sum_{n=0}^{\infty}  \frac{\tau^{2n+1}}{(2n+1)!} C^n.
\end{displaymath}
Using only $m$ products, $C^2, C^3,\ldots, C^{m+1}$, we can compute the truncation
\begin{equation*}
E_2^{[m]}=\begin{pmatrix} 
\sigma_m & 
\mu_m \\ 
-C\mu_m & 
\sigma_m   \end{pmatrix}
\text{ with }
\sigma_m =\sum_{n=0}^{m+1} \frac{\tau^{2n}}{(2n)!} C^n, \;
\mu_m =\sum_{n=0}^{m}  \frac{\tau^{2n+1}}{(2n+1)!} C^n,
\end{equation*}
and commit an error of $E_2-E_2^{[m]}=\mathcal{O}(\tau^{2m+3})$.
The result is not a symplectic matrix, but can be made symplectic by adding a correction term $\delta_m$ as follows
\begin{equation}\label{approxe3}
E_2^{[m,s]}=\begin{pmatrix} 
\sigma_m & 
\mu_m +\delta_m \\ 
\nu_m & 
\sigma_m   \end{pmatrix}
\end{equation}
with $\delta_m$ a close-to-zero symmetric matrix. The symplecticity condition for $E_2^{[m,s]}$ reads
\[
  (E_2^{[m,s]})^T\,J\,E_2^{[m,s]}=J.
\]
Since $C$ is a symmetric matrix, $\sigma_m$ and $\mu_m$ are symmetric and commute, $\sigma_m\mu_m-\mu_m\sigma_m$, and the symplectic condition simplifies to 
\[
  \sigma_m^2 - (\mu_m +\delta_m)\nu_m = I 
	\quad \Longrightarrow \quad
	\delta_m= (\sigma_m^2 - I)/\nu_m - \mu_m.
\] 
Hence, $\delta_m$ can be computed with a product and one inverse for a total extra cost of $(1+\frac43) \CC$. Summarizing, we can achieve a symplectic approximation of error $\mathcal{O}(\tau^{2m+1})$ with $m+\frac43$ products.

%

Notice that the local truncation error is then
\[
 {\cal O}(\tau^{2m+1})=\gamma^{2m+1}{\cal O}(h^{2m+1})
\]
and methods with small values of the coefficients $\gamma$ in the composition can approximate the exponentials using a lower order truncation and hence at lower computational cost.

\subsection*{Acknowledgments}

PB and SB acknowledge the Ministerio de Econom\'{\i}a y Competitividad
(Spain) for financial support through the coordinated project MTM2013-46553-C3.

\end{document}